\documentclass[12pt, a4paper]{article}
\usepackage[margin=25.4mm]{geometry}
\usepackage[T1]{fontenc}
\usepackage{newtxtext, newtxmath}
\usepackage[ruled]{algorithm2e}
\usepackage{titlesec}
\titleformat*{\section}{\raggedright\large\bfseries}
\titleformat*{\subsection}{\raggedright\normalsize\bfseries}
\usepackage[colorlinks]{hyperref}
\usepackage[numbers, sort&compress]{natbib}
\usepackage{paralist}
\usepackage[normalem]{ulem}
\renewcommand*{\underline}[2]{\uline{#1} (Origin: #2)}
\usepackage[english]{babel}                     % Language setting
\usepackage{blindtext}
\usepackage{setspace}
\usepackage{lscape}
\usepackage{array}
\usepackage{hyperref}
\usepackage{appendix}
\usepackage{autobreak}
\hypersetup{hidelinks} 
\usepackage{ breqn}
\usepackage{autobreak}
\usepackage{graphicx}                           % Graphics
\usepackage{amsmath}                            % Math
\usepackage{hyperref}
\usepackage{algorithmicx}
\usepackage{lastpage}                           % Footer note
\usepackage{lipsum}                             % For dummy text
\usepackage{indentfirst}
\usepackage{amsfonts} %\mathbb
\usepackage{amsopn}
\usepackage{array}
\allowdisplaybreaks
\begin{document}
	\begin{flushleft}\Large\textbf{
			Four Limit Cycles in Three-Dimensional Competitive Lotka-Volterra  Systems of Class 28 in Zeeman's Classification}
		
	\end{flushleft}
	\begin{flushleft}
		\textbf{Mingzhi Hu$^{1}$, Zhengyi Lu$^{1}$ and Yong Luo$^{2}$}\\[1em]
		
		%\begin{compactenum}[$^1$]
		$^{1}$ School of Mathematical Sciences, Sichuan Normal University,  Chengdu 610068, China
		
		$^{2}$ College of Mathematics and Physics, Wenzhou University,  Wenzhou 325000, China
		%\end{compactenum}
		
		\hspace{1em}
		
		Correspondence should be addressed to Yong Luo; luoyong@wzu.edu.cn
	\end{flushleft}
	
	\section*{Abstract}
	In this paper, a three-dimensional competitive Lotka-Volterra system with four limit cycles is constructed for class 28 in Zeeman's classification. Combined with existing results --- from Gyllenberg and Yan (2009) for class 27, from Wang, Huang and Wu (2011) for class 29, and from Yu, Han and Xiao (2016) for class 26 --- our finding indicates that there exist systems with at least four limit cycles for each class among classes 26 $-$ 29. %This work provides a partial answer to the conjecture proposed in Hofbauer and So (1994) and further extended in Yu, Han and Xiao (2016).}%together with the results from Gyllenberg, Yan and Wang (2009) for class 27, from Wang, Huang and Wu (2011) for classes 28 and 29 and Yu, Han and Xiao (2016) for class 26 which indicate that for each class among classes 26$-$29, there exist systems with  at least four limit cycles. This gives a partial answer to a conjecture  proposed in  Hofbauer and So (1994) as well as in Yu, Han and Xiao (2016).}      % the abstract
	\vspace{1em}

\noindent Keywords: Lotka-Volterra; competitive systems; Zeeman's classification; class 28; four limit cycles   % the 
\vspace{1em}

\section{Introduction}

MacArthur \cite{Mc} constructed a Lyapunov function for a  Lotka-Volterra system with a symmetric interaction matrix, complementing the original Lyapunov function proposed by Volterra \cite{Vol}. By means of MacArthur’s Lyapunov function, it can be shown that an $\omega$-limit set for a competitive system with a symmetric interaction matrix is either an isolated equilibrium or a continuum of equilibria.
%Using MacArthur's Lyapunov function, it can be shown that an $\omega-$limit set of a competitive system with a symmetric interaction matrix is either an isolated equilibrium or a continuum of equilibria.

In 1975, May and Leonard \cite{May} presented a landmark result for three-dimensional competitive systems. They constructed a three-dimensional competitive Lotka-Volterra system with a cyclically symmetric interaction matrix and obtained a repulsive heteroclinic cycle.%the periodicity of the system under certain conditions. 
Meanwhile, through numerical simulation, Gilpin \cite{Gil} found that a  three-dimensional competitive system may admit a limit cycle.
%In 1975, May and Leonard \cite{May} have achieved a pioneering landmark result for three-dimensional competitive Lotka-Volterra systems. They constructed a cyclically symmetric three-dimensional competitive Lotka-Volterra system with uniform intrinsic growth rates and symmetric competition coefficients, and rigorously derived the periodic dynamical behavior of the system under the specific parameter condition of $\alpha+\beta=2$. Meanwhile, via numerical simulation, they first identified that such a cyclically symmetric three-dimensional competitive system could admit a stable periodic limit cycle and an attractive or repulsive heteroclinic cycle, laying a foundational framework for subsequent studies on limit cycles in high-dimensional competitive Lotka-Volterra systems. In the same year, Gilpin \cite{Gil} independently conducted numerical investigations on limit cycle phenomena in competition communities with three or more competing species, with a research focus on the biological implications of such cyclic dynamical behaviors, thus forming an important complement to the mathematically oriented analysis of May and Leonard \cite{May}.

In 1979,  Coste, Peyraud and Coullet \cite{CPC} and Schuster, Sigmund and Wolff \cite{SSR} independently proved the existence of a limit cycle for three-dimensional competitive systems (which  actually belong to class 27 in Zeeman's classification \cite{Zee}), and verified the attractiveness or repulsion of the corresponding heteroclinic cycle. In 1981, Hofbauer \cite{Ho81} proved the existence of limit cycles in general $n$-dimensional $(n \geqslant 4)$ competitive systems via Hopf bifurcation theory.

Smale \cite{Smale} in  1976  and  Hirsch \cite{Hir88} in 1988 demonstrated that the dynamical behaviors of $n$-dimensional competitive systems are equivalent to those of  $(n-1)$-dimensional systems.%showed that the dynamical behavior of $n-$ dimensional competitive systems is comparable to that of $(n-1)-$ dimensional systems. 
For a three-dimensional competitive Lotka-Volterra  system, Hirsch's result \cite{Hir88} ensures that there exists an invariant manifold (called a carrying simplex), which is homeomorphic to a two-dimensional simplex and attracts all orbits except the origin. Based on Hirsch's theorem, Zeeman \cite{Zee} used the qualitative analysis of the system at the boundary surface to define the combinatorial equivalence relation through parameter inequalities, and obtained thirty-three stable classes of  three-dimensional Lotka-Volterra competitive systems.  Zeeman's result indicates that the limit sets for 27 classes among the 33 ones are all fixed points \cite{Zee}, so the dynamical behaviors of systems in these 27 classes are fully described. Using the Hopf bifurcation theorem, Zeeman further proved that for each of the remaining six classes (class 26 to class 31), there always exists a system that can be constructed to admit an isolated periodic orbit or limit cycles via appropriate parameter selection.%Zeeman further proved that among each remaining the six classes: from class 26 to class 31, there is always a system that can be constructed to have an isolated periodic orbit or limit cycles by selecting parameters.

The first result of a three-dimensional competitive Lotka-Volterra  system with two limit cycles was obtained in class 27 of Zeeman's classification by Hofbauer and So \cite{Ho} based on Hirsch's monotonic flow theorem, the center manifold theorem and the Hopf bifurcation theorem. In their example, a locally stable positive equilibrium is surrounded by two limit cycles: one is from the Hopf bifurcation theorem, and the other is guaranteed via the Poincaré-Bendixson theorem. In \cite{Ho}, Hofbauer and So proposed a problem of how many limit cycles may exist in the six classes 26$-$31 in Zeeman's classification.

In 2002, Lu and Luo \cite{Lu02} constructed two limit cycles in class 26, class 28, and class 29, respectively, which partially answered Hofbauer and So's question about the existence of two limit cycles in other classes.

In 2003, the first attempt to construct three limit cycles was made by Lu and Luo \cite{Lu03}, who proposed a system with three limit cycles in Zeeman's class 27. However, as Yu, Han and Xiao \cite{YP} pointed out, due to a sign error, the positive definiteness of the Lyapunov function for calculating the focal values was not guaranteed, and thus Lu and Luo's original attempt was invalid.

%In 2003, the first attempt to construct three limit cycles was given by Lu and Luo \cite{Lu03}. They constructed three limit cycles in Zeeman's system of 27 classes. As Yu, Han and Xiao \cite{YP} pointed out that due to a sign error, the positive definiteness of the Lyapunov function for calculating the focal value was not  guaranteed, therefore, Lu and Luo's attempt failed.

In 2006, Gyllenberg, Yan and Wang \cite{GY2006} found that the system with two limit cycles in class 29 constructed by Lu and Luo \cite{Lu02} can actually admit a third limit cycle. The instability of the outer small-amplitude limit cycle, together with the repulsion of the boundary of the carrying simplex, satisfies the conditions of the Poincaré-Bendixson theorem, which guarantees the existence of a third limit cycle.

%In 2008, Lian, Lu and Luo \cite{lian} found through automatic search that both class 30 and class 31 in Zeeman's classification can have systems with three limit cycles. Two of them are small-amplitude limit cycles, and the other is a large-scale limit cycle obtained via the Poincaré-Bendixson theorem. %Recently, Li and Jiang \cite{jiang}  pointed out their conclusion about stability was wrong, actually, it seems that there was a typing error since three limit cycles can be constructed by another answer. In fact, for $\lambda\in[-\frac{1919094527}{34359738368},-\frac{144294119}{2583462210}]$ in their example, we can obtain an unstable 2-order focal value which ensures the existence of two limit cycles. Finally, a third limit cycle can be obtained by applying Poincaré-Bendixson theorem.

There were also works by Gyllenberg and Yan \cite{GY2007} in 2009 and Wang, Huang and Wu \cite{w} in 2011, respectively, for the existence of multiple limit cycles for class 30. Gyllenberg and Yan claimed to have constructed two limit cycles in class 30. One is a small-amplitude limit cycle obtained via Hopf bifurcation, and the other is obtained via the Poincaré-Bendixson theorem. As pointed out by Yu, Han and Xiao \cite{YP}, their construction of limit cycles was incorrect since the Lyapunov function they gave is not positive definite. However, with careful calculation, three limit cycles can actually be constructed by using their original  example \cite{GY2007} \cite{YP}. First, two small-amplitude limit cycles can be constructed such that the outer one is stable. Since the system belongs to class 30, the boundary of the carrying simplex is an attractor; therefore, the third limit cycle can be obtained via the Poincaré-Bendixson theorem.  Three years later, Wang, Huang, and Wu \cite{w} proved that for each class in Zeeman's classes 26$-$31, there exists a system with at least three small-amplitude limit cycles. Unfortunately, their result for class 30 was invalid, as the constructed system was not competitive.

In 2009, Gyllenberg and Yan \cite{GY2009} constructed an example in class 27 similar to \cite{Lu03} and claimed that they obtained four limit cycles, three of which were obtained due to the Hopf bifurcation, and the fourth one was obtained through the existence of the heteroclinic cycle and using the Poincaré-Bendixson theorem. Similarly, Gyllenberg and Yan neglected the positive definiteness of the Lyapunov function. Yu, Han and Xiao \cite{YP} pointed out that there does exist an example in \cite{GY2009} belonging to class 27, which can admit four limit cycles. A system in class 29 with four limit cycles was given by Wang, Huang, and Wu \cite{w}.

Recently, Yu, Han, and Xiao \cite{YP} gave two examples  in class 27 and two in class 26 with four small-amplitude  limit cycles, respectively. Besides, Li and Jiang \cite{jiang} showed three limit cycles for classes 28 and 31.

In summary, multiple limit cycles can appear in each of the six classes of  Zeeman's classification: classes 26, 27, 28, 29, 30 and 31. The known results to date are as follows: for the classes 26, 27 and 29, there are at least four limit cycles.

An open problem proposed by Yu, Han and Xiao \cite{YP} is whether there are four limit cycles in all classes 26---31.

In this paper, by combining the limit cycle construction algorithm proposed by Hofbauer and So that modified by Lu and Luo \cite{Lu02} together with the real root isolation algorithm presented in \cite{shigen}, we construct a system with four limit cycles in class 28.%30 and 31 in Theorem 3.1 and Theorem 3.2 of the present paper by automatic search. By using the program \texttt{3DLVzd} written by the authors, 6733 examples are searched with randomly chosen interaction matrices. %And among these examples, 350, 3, 21, 1 and 1 limit cycles belong to classes 27, 28, 29, 30 and 31, respectively. This gives an affirmative answer to Yu, Han and Xiao's problem \cite{YP}.

Summarizing the results, we have the following main result.

\noindent\bf{Theorem 1.1.} \normalfont\emph{For each of the four classes 26, 27, 28 and 29 in Zeeman's classification, there exists a system with four limit cycles.}
%A series of the attempts, challenges, mistakes and corrections in the previous and current work show the following conclusion.

%Theorem: For each class of the six ones  (Class 26, Class 27, Class 28, Class 29, Class 30, and Class 31)  in Zeeman classification,  there exists a system with 4 limit cycles.

%\newpage
\section{Automated search algorithm}
In this section, we adopt the algorithmic construction method proposed by Hofbauer and So \cite{Ho} and modified by Lu and Luo \cite{Lu02} to search for limit cycles.

Consider a three-dimensional Lotka-Volterra system

\begin{equation}
\begin{aligned}
	\dot{x}_i=x_i(\sum_{j=1}^{3}a_{ij}(x_j-1)),
\end{aligned}
\end{equation}
where $a_{ij} <0$ and $\textbf{1}= (1,1,1)$ is the unique positive equilibrium of system (1). %System (1) can be expressed in vector form as follows
%\begin{equation}
%	\begin{aligned}
%		\dot{x}=diag(x)A(x-\mathbf{1}).
%	\end{aligned}
%\end{equation}

%where $x= (x_1, x_2, x_3), A={(a_{ij}}})_{33}, x-1=(x_1-1, x_2-1, x_3-1).$

Suppose that the matrix $A = (a_{ij})_{3\times3}$ has a real eigenvalue $\lambda $ and a pair of purely imaginary eigenvalues $\pm \omega i$ $(\omega\ne 0)$. Then there exists a transformation matrix $T$ that transforms $A$ into the block diagonal form. $$TAT^{-1}=\left[\begin{array}{ccc}
c_{11} & c_{12} & 0 
\\
c_{21} & c_{22} & 0 
\\
0 & 0 & \lambda  
\end{array}\right].
$$
Here, the submatrix $\left(\begin{array}{cc}
c_{11} & c_{12} 
\\
c_{21} & c_{22} 
\end{array}\right)
$ has a pair of purely imaginary eigenvalues $\pm\omega i$$(\omega\ne0)$; that is, the submatrix satisfies $c_{11}+c_{22}= 0 $ and $c_{11}c_{22}-c_{12}c_{21}>0$.
%In addition, to guarantee the positive definiteness of the Lyapunov function, we require that $c_{21}<0$.

From the center manifold theorem \cite{Carr}, we can suppose that, under the transformation $y=T(x-1)$, the transformed system with linear part $Cy$ has an approximation to the center manifold taking the form

$y_3 = h(y_1,y_2) = h_2(y_1, y_2) + h_3(y_1, y_2) + h_4(y_1, y_2) + h_5(y_1, y_2) + h_6(y_1, y_2) + \text{h.o.t},$\\
where $ y=(y_1,y_2,y_3)^T, h_i(y_1,y_2)=\sum_{j=0}\limits^{i}p_{ij}y_1^{j}y_2^{i-j},(i=2,\dots,6)$ and %$y_1$ and $y_2$ are homogeneous polynomial of degree k. 
h.o.t. denotes the terms with orders greater than or equal to seven.

We have the following steps to calculate the focal value to get the relevant conclusion:

\noindent \begin{tabular}{lp{12.5cm}}
Step 1:  & Using the subprograms \texttt{randA} and \texttt{DIA (A = randA(), T = DIA(A))}, the coefficient matrix $A$ is  randomly selected, and the block diagonal matrix $TAT^{-1}$ is obtained. \\
Step 2: &The coefficients $p_{ij}$ are solved via the subprogram \texttt{Vh}, and thus the approximate center manifold $y_3 = h(y_1, y_2)$ is obtained. \\

Step 3: &By substituting $y_3 = h(y_1, y_2)$ into $\dot{y}_1$ and $\dot{y}_2$, a two-dimensional system of center-focus type is obtained. \\

Step 4: &The focal values are calculated via the subprogram \texttt{JJLLineSolve}. \\

Step 5: &The real root isolation algorithm (subprogram \texttt{mrealroot} \cite{shigen}) is applied to verify the independence of the focal values and construct multiple limit cycles.\\

Step 6: &The subprogram \texttt{CL} is used to determine which class in Zeeman's classification the constructed system belongs to.\\
\end{tabular}
%\vspace{1em}
\newpage
The main algorithm of the program \texttt{3DLVzd} is as follows:
\begin{algorithm}[!ht]
\KwIn{number}
\KwOut{coefficient matrix A, Zclass}
\Begin{
\For {$zz$ $to$ $number$}	
{
	{restart:with(linalg):with(LinearAlgebra):\\}
	{A:=randA(A):P :=DIA(A, $\mu$, [y[1], y[2], y[3]]):\\}
	{g1 := map(factor, P[1]):g2 := map(factor, P[2]):g3 := map(factor, P[3]):\\}
	{hh := factor(Vh2(g1, g2, g3, [y[1], y[2], y[3]], h, 6)):\\}
	{f1 := subs(y[3] = hh, g1):f2 := subs(y[3] = hh, g2):\\}
	{La := factor(JJLLinearSolve(f1, f2, [y[1], y[2]], 3)) :\\}
	{la(i) := expand(numer(La[i])):(i=1,2,3):lad(i) := expand(denom(La[i])):(i=1,2,3):\\}
	{sys :=[la1, la2]:X := [$\lambda$,n]:\\}
	{S := charsets[mcharset](sys, X, basset):\\}
	{S := S[1]:charsets[iniset](S, X):\\}
	{S1 := subs($\lambda$ = -l, S[1][1]):\\}
	{SS := factor(subs(n = -q, S[1][2])):\\}
	{s1 := realroot(S1,$1/10^{10}$):\\}
	{J:=mrealroot(S,{$\lambda$,n}):Check:=CC(A,o,o1):\\}
	{B := subs({n = o,$\lambda$ = o1,$\mu$ = o2}, A):P1 := CL(B):\\}
	{\If{P1[1]=Zclass}
		{\If{$LV_3\in$ number and Check==true}
			{save A,P1[1],"zd.txt":}}}
}}
\caption{\texttt{Automated search algorithm}.\label{alg:1}}
\end{algorithm}

%Compared with \cite{lian}, the number of terms and degree in each of the focal values after factorization are given below:

%\begin{center}
%\begin{tabular}{c|c|c}
%\hline
%Focal Value &\cite{lian}(Number of %Terms,Degree)&(Number of Terms,Degree)\\\hline
%$LV_1$&(3,4)&(3,8)\\
%$LV_2$&(14,13)&(63,26)\\
%$LV_3$&None&(169,50)\\\hline
%\end{tabular}\end{center}
Compared with \cite{Lu02}, we observe that the number of terms and the degree of the multivariate polynomials increase very rapidly as the order of the focal values increases.

On the other hand, we focus on much more complicated multivariate polynomial systems in this paper, based on the real root isolation algorithm for multivariate polynomial systems presented in \cite{2005}. In our program, a system of multivariate polynomials is triangularized by using the Epsilon package given by Wang \cite{WWD}. 

\section{Four limit cycles in class 28}
The dynamical behavior of systems of the class 28 in Zeeman's classification restricted to the carrying simplex is as follows:
\begin{center}
\includegraphics[scale=0.1]{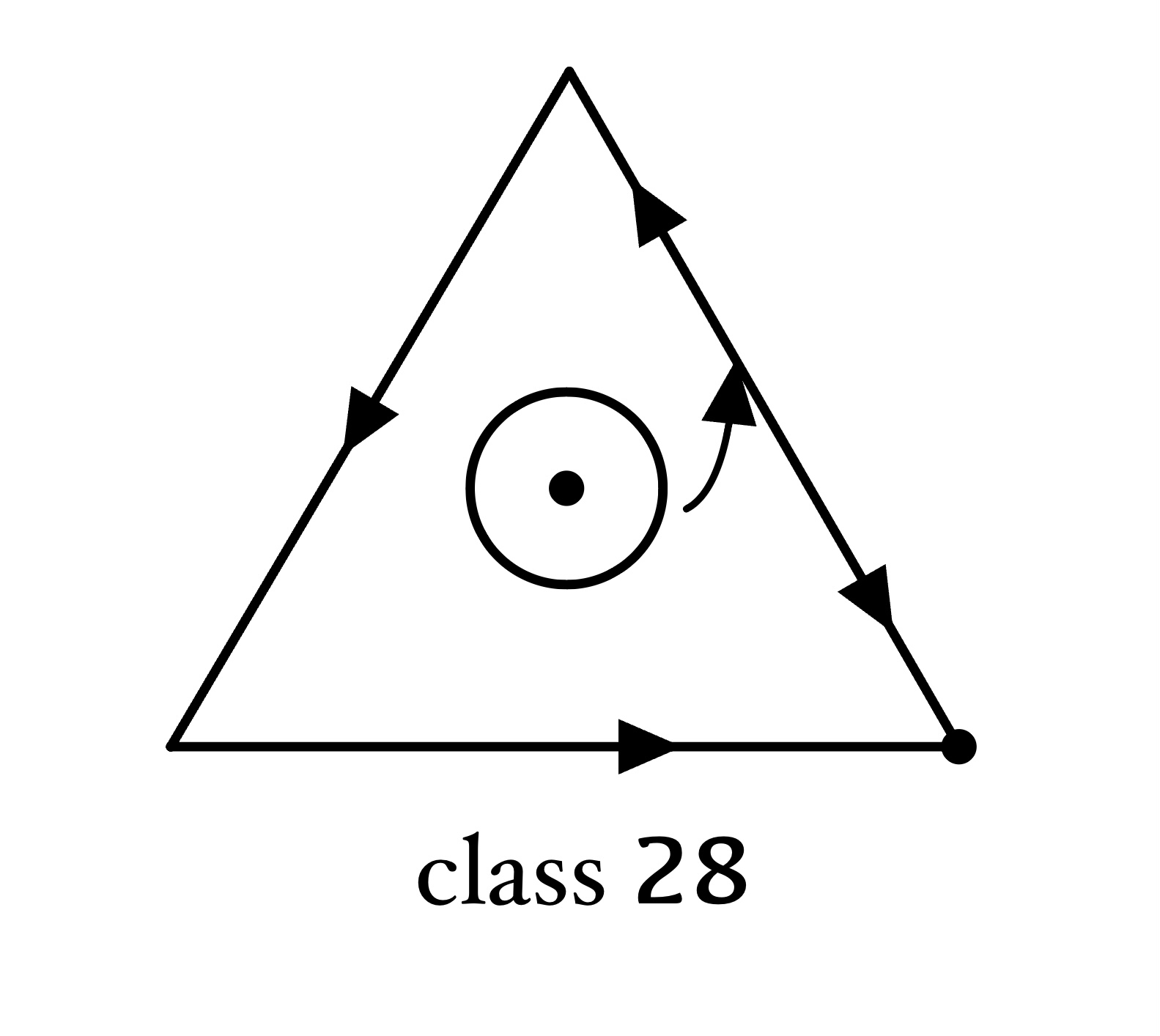}
\end{center}

The interaction matrix of the constructed system is chosen by:

$$A=\left[\begin{array}{ccc}
-\frac{17}{24} & -2 & -\lambda  
\\
-\frac{33}{23} & -10 & -\mu  
\\
-n  & -21 n  & -\frac{99}{37} 
\end{array}\right],(n,\mu,\lambda>0).$$

To satisfy the necessary eigenvalue condition \cite{Ho} $\text{det}(A) = (A_{11} + A_{22} + A_{33}) \cdot \text{trace}(A)$, we need $\mu=-\frac{607835112 \lambda  n -7773334823}{4864016448 n}$. 

By the transformation $y = T (x-1) $ with $T$ as follows,\begin{spacing}{2}$$\left[\begin{array}{ccc}
-\frac{33}{23} & \frac{3005}{888} & \frac{607835112 \lambda  n -7773334823}{4864016448 n} 
\\
-n  & -21 n  & \frac{257}{24} 
\\
\frac{175916243}{65729952}+\frac{4791493 \lambda  n}{1825832} & 21 \lambda  n +\frac{257}{12} & \frac{7622102628 \lambda  n +7773334823}{2432008224 n} 
\end{array}\right].
$$
\end{spacing} 
The three-dimensional system is transformed into a new one whose linear part is in the block diagonal form
\begin{center}\begin{spacing}{2}$\left[\begin{array}{ccc}
	\frac{1014026}{684499} & \frac{416062526328888 \lambda  n +384134040047899}{3329414394639552 n} & 0 
	\\
	-\frac{8896983 n}{684499} & -\frac{1014026}{684499} & 0 
	\\
	0 & 0 & -\frac{11885}{888} 
\end{array}\right]
$.\end{spacing}\end{center}
Furthermore, the three-dimensional system can be reduced to a two-dimensional system by the subprogram $Vh$, and
the first three focal values can be obtained via the subprogram \texttt{JJLLineSolve},
$$\begin{aligned}LV_1=f(\lambda,n)=\frac{f_1(\lambda,n)}{f_2(\lambda,n)},\\
LV_2=g(\lambda,n)=\frac{g_1(\lambda,n)}{g_2(\lambda,n)},\\
LV_3=z(\lambda,n)=\frac{z_1(\lambda,n)}{z_2(\lambda,n)},\end{aligned}$$
where
\begin{align*}
f_1(\lambda,n) = &-684499 (71175864 n \lambda -30452821)\\& (1867427763509790559220459722237440 \lambda^{3} n^{3}\\&-8066558192490463098597559057769472 \lambda^{2} n^{3}\\&-9284198345879360722318571367478464 \lambda^{2} n^{2}\\&+9630923381872490306204845292994048 \lambda \,n^{3}\\&+18212185214244398672238510809517312 \lambda \,n^{2}\\&+5385169712442285368618043473601672 n \lambda \\&+37509227186769280161709353461815488 n^{2}\\&+1036086857152915319628573370644784 n \\&-604462354449619944145534311192809),\\\\

f_2(\lambda,n)=&7401027927168 n^{3} (1169277093792 \lambda  n +31737803345137)\\& (146159636724 \lambda  n +16056506276339)^{2},\\\\

g_2(\lambda,n)=&2717701945572392898631507132761600 (1315436730516 \lambda  n +15556227332951)\\& (1169277093792 \lambda  n +31737803345137)^{3} (146159636724 \lambda  n +16056506276339)^{6} \\&(416062526328888 \lambda  n +384134040047899)^{2} n^{5}

\\\\

z_2(\lambda,n)=&74656592898569603821499140755938066449603169746217765376000 \\&(4677108375168 \lambda  n +30236966514973) (146159636724 \lambda  n +16056506276339)^{10}\\& (1169277093792 \lambda  n +31737803345137)^{5} (1315436730516 \lambda  n +15556227332951)^{2} \\&(71175864 \lambda  n -30452821) (416062526328888 \lambda  n +384134040047899)^{5} n^{7}
\end{align*}

%Here, the expressions of $g_1(\lambda,n)$ (a polynomial of $63$ terms with degree $26$) and $z_1(\lambda,n)$ (a polynomial of $169$ terms with degree $50$) are listed in the appendix.
Here, $g_1(\lambda,n)$ is a polynomial of 63 terms with degree 26, and $z_1(\lambda,n)$ is a polynomial of 169 terms with degree 50.

By running the program \texttt{mrealroot} \cite{shigen}, we obtain:

\begin{center}
>-\texttt{mrealroot}$([f_1(\lambda,n),g_1(\lambda,n)],[\lambda,n],\displaystyle\frac {1}{10^{20}},[z_1(\lambda,n),z_2(\lambda,n),$\text{det}(A)]),
\end{center}
with

\begin{small}
%$\begin{spacing}{1.5}
$$[[\lambda_1,n_1][+,-,-]],$$
%$$[[\lambda_2,n_2][+,-,-]],$$
such that
$$\lambda_1\in({\frac{48083713211257141381227}{9444732965739290427392}}, {\frac
	{48083713211499877152963}{9444732965739290427392}}),$$ 
and
$$n_1\in({\frac{
		18617876387518095278417715070016705816645420386546066217507077513538675027145
	}{
		57896044618658097711785492504343953926634992332820282019728792003956564819968
}},$$

$$ {\frac{
		9308938193759047639208857535008352908322710193273033108753538756769337513573
	}{
		28948022309329048855892746252171976963317496166410141009864396001978282409984
}}).$$

%\end{spacing}
\end{small}
\bigskip

This result confirms that the constructed system is competitive, and the two real roots of $LV_1 =LV_2=0$ in the interval form leads that satisfies $LV_3<0.$

Thus, we can perturb $LV_1$ and $LV_2$ by using $n$ and $\lambda$ to obtain two small-amplitude limit cycles. Further, we change $\mu$ to perturb $LV_0$ such that $LV_0\cdot LV_1 < 0$ and $|LV_0| \ll |LV_1|$ to obtain an additional small-amplitude limit cycle. This ensures the existence of three small-amplitude limit cycles, and the outermost one is stable.

Finally, we verify whether the Poincaré-Bendixson theorem can be applied to obtain a fourth limit cycle. Using Zeeman's notation \cite{Zee}, we have $R_{ij} =\text{sgn}(\alpha_{ij})$  and $Q_{kk}=\textrm{sgn}(\beta_{kk})$, with $\alpha_{ij}=\frac{b_ia_{ji}}{a_{ii}}-b_j=(AR_i)_j-b_j$ and $\beta_{kk}=(AQ_k)_k-b_k$, which are the algebraic invariants of $A$. Here, $R_i$ is the equilibrium on the $x_i$-axis, and $Q_k$ is the positive equilibrium on the plane $x_k=0$.

Since$$R_{12}=Q_{33}=R_{21}=-R_{23}=R_{32}=-R_{31}=R_{13}=1,$$this implies that the constructed example belongs to class 28 in Zeeman's classification. Since the third focal value is negative, the outer limit cycle is stable. For class 28, the boundary of the carrying simplex is an attractor, which guarantees the existence of a fourth limit cycle by the Poincaré-Bendixson theorem. Hence, we have

\noindent\bf{Theorem 3.1.} \normalfont\emph{There exist at least four limit cycles for class 28 in Zeeman's classification.}

\section{Concluding remarks}
In this paper, three-dimensional competitive Lotka-Volterra systems with four limit cycles for Zeeman's class 28 are constructed. By incorporating Hirsch's dimension reduction method, the center manifold construction, the focal values computation   combined with the real root isolation algorithm  with multivariable polynomials \cite{2001,2005} and automated search, we demonstrate the independence of the focal values and automatically verify the stability consistency between the outermost small-amplitude limit cycle and the boundary of the carrying simplex. Thus, three small-amplitude limit cycles  and a large-scale limit cycle are obtained. Therefore, for systems in Zeeman's class 28, the existence of four limit cycles is confirmed.

It seems that to construct four limit cycles for systems in Zeeman's classes 30 and 31 is an interesting and challenging problem. After deriving focal values up to the third order, the remaining challenge is to triangulate these high-degree polynomials to verify their independence and ensure stability consistency between the outermost small-amplitude limit cycle and the boundary of the carrying simplex. Directly calculating focal values up to the fourth order and proving their independence may be more complicated for these two cases. We will leave these problems for future consideration.

%\newpage
\end{document}